\documentclass[12pt]{amsart}
\font\emailfont=cmtt10

\usepackage{amsmath,amsthm,amsfonts,amscd,flafter,epsf}
\usepackage{graphics}
\usepackage{epsfig}
\usepackage{psfrag}

\newcommand\commentable[1]{#1}

\newtheorem{theorem}{Theorem}[section]
\newtheorem{prop}[theorem]{Proposition}
\newtheorem{cor}[theorem]{Corollary}

\def\endproof{\relax\ifmmode\expandafter\endproofmath\else
  \unskip\nobreak\hfil\penalty50\hskip.75em\hbox{}\nobreak\hfil\bull
  {\parfillskip=0pt \finalhyphendemerits=0 \bigbreak}\fi}
\def\endproofmath$${\eqno\bull$$\bigbreak}
\def\bull{\vbox{\hrule\hbox{\vrule\kern3pt\vbox{\kern6pt}\kern3pt\vrule}\hrule}}

\newcommand{\C}{\mathbb{C}}

\newcommand{\Z}{\mathbb{Z}}

\newcommand{\cm}{\cdot}

\newcommand\relspinc{\underline{\spinc}}

\newcommand\Filt{\mathcal F}

\newcommand\ModSphere{\ModFlow\left({\mathbb S}\longrightarrow 
\Sym^{g-1}(\Sigma_{1})\times \Sym^2(\Sigma_{2})\right)}
\newcommand\ModSpheres\ModSphere

\newcommand\CFa{\widehat{CF}}

\newcommand\HFa{\widehat{HF}}

\newcommand\UnparModSp{\widehat \ModSp}
\newcommand\UnparModFlow\UnparModSp
\newcommand\Mod\ModSp

\newcommand{\spinc}{\mathfrak s}

\newcommand\ModMaps{\mathcal M}
\newcommand\ModSp\ModMaps

\newcommand\spincrel\relspinc

\newcommand\Dual{\mathcal D}
\newcommand\Duality\Dual

\newcommand\ons{Ozsv{\'a}th and Szab{\'o}}
\newcommand\os{Ozsv{\'a}th-Szab{\'o}}

\commentable{

\title[{Notions of positivity and The Ozsv{\'a}th-Szab{\'o} concordance invariant}] 
{Notions of positivity and the Ozsv{\'a}th-Szab{\'o} concordance invariant}

}

\author[Matthew Hedden]{Matthew Hedden}
\address{Department of
Mathematics, Princeton University, NJ \newline
\indent{\emailfont{mhedden@math.princeton.edu}}}

\begin{document}

\begin{abstract}

In this paper we examine the relationship between various types of positivity for knots and the concodance invariant $\tau$ discovered by \ons \ and independently by Rasmussen.  The main result shows that, for fibered knots, $\tau$ characterizes strong quasipositivity.  This is quantified by the statement that for $K$ fibered, $\tau(K)=g(K)$ if and only if $K$ is strongly quasipositive.   In addition, we survey existing results regarding $\tau$ and forms of positivity and highlight several consequences concerning the types of knots which are (strongly) (quasi) positive.  For instance, we show that any knot known to admit a lens space surgery is strongly quasipositive and exhibit infinite families of knots which are not quasipositive.
	
\end{abstract}

\maketitle
\section{Introduction and Background}

There are many notions of positivity for braids and knots.  Perhaps simplest is that of a {\em positive} knot.  A knot is said to be positive if it has a projection for which the writhe equals the crossing number \cite{Lickorish}.  A weaker notion is that of a {\em positive braid}.  Positive braids are those knots and links which can be obtained as the closure of a word in the braid group consisting only of positive generators, $\Pi_{k=1}^m\sigma_{i_k}$ (i.e. there are no $\sigma_i^{-1}$).  Of course a positive braid is positive as a knot.  Perhaps slightly more abstruse are the concepts of quasipositivity and strong quasipositivity.  These concepts have been studied extensively by Rudolph and the present paper draws heavily on his collected works.  For a thorough introduction to the subject see \cite{Rudolph1}.  Our motivation is to better understand the relationship between different notions of positivity - particulary strong quasipositivity - and the concordance invariant $\tau$ coming from \os \ homology.  This invariant was discovered by \ons \ and independently by Rasmussen in his thesis, see \cite{FourBall,Ras1}. (This invariant should not be confused with Rasmussen's $s$ invariant \cite{Ras2} - see remarks belo).

\begin{figure}

\begin{center}
	\psfrag{a}{$S_{2,5}$}
	\psfrag{b}{$S_{1,5}$}
	\psfrag{c}{$S_{2,4}$}
	\psfrag{d}{$S_{1,2}$}

\includegraphics{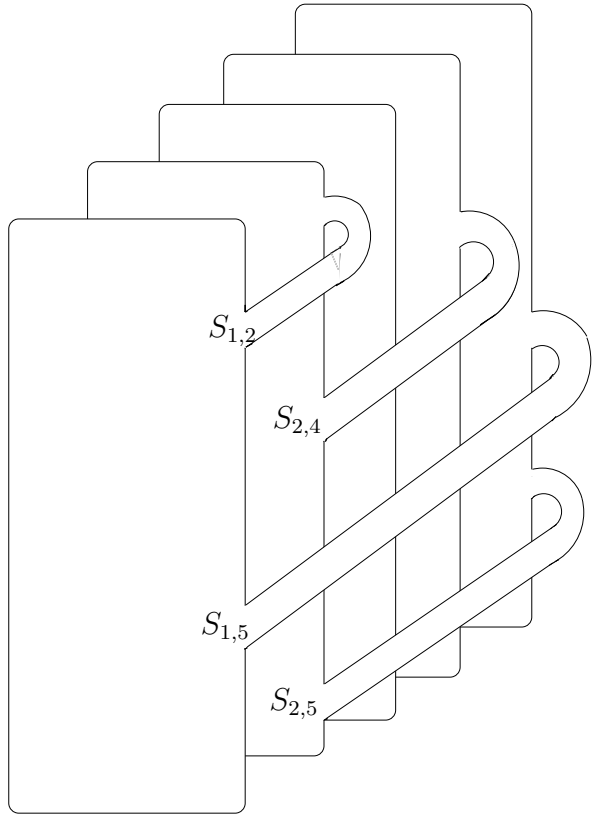}

\end{center}

\caption{\label{fig:SS} Quasipositive Seifert surface.  The boundary of this surface is the strongly quasipositive knot realized as the closure of the braid $\beta=\sigma_{2,5}\sigma_{1,5}\sigma_{2,4}\sigma_{1,2}$.  It is obtained from $5$ parallel disks by attaching $4$ positive bands, $S_{i,j}$.  Each band starts from the left and curves up as it attaches at the right.  A negative band (not shown) would curve down at the right. It should be clear that attaching the band $S_{i,j}$ is equivalent to adding $\sigma_{i,j}$ to the braid $\beta$. }

\end{figure}

A {\em strongly quasipositive} knot is a knot which has a special kind of Seifert surface, a so-called {\em quasipositive Seifert surface}.  Quasipositive Seifert surfaces are those surfaces obtained from $n$ parallel disks by attaching positive bands.  The easiest way to define positive bands and quasipositive Seifert surfaces is through figures and so we refer the reader to Figure \ref{fig:SS} for an illustration.  A knot or link is strongly quasipositive if it can be realized as the boundary of such a surface. \newpage Let $B_n$ denoted the braid group on $n$ strands, with generators $\sigma_1,\ldots,\sigma_{n-1}$, and let

$$\sigma_{i,j}=(\sigma_i\ldots\sigma_{j-2})(\sigma_{j-1})(\sigma_i\ldots\sigma_{j-2})^{-1}.$$

It is evident from Figure \ref{fig:SS} that strongly quasipositive knots are precisely those knots which can be realized as the closure of braids of the following form:

$$\beta = \Pi_{k=1}^m \sigma_{i_k,j_k}.$$

Noting that $\sigma_{i,j}$ is of the form $w\sigma_{j-1}w^{-1}$ with $w\in B_n$, the weaker notion of a {\em quasipositive}  knot is any knot which can be realized as the closure of a braid of the form:

$$\beta = \Pi_{k=1}^m w_k \sigma_{i_k} w_k^{-1}.$$

Thus quasipositive knots are closures of braids consisting of arbitrary conjugates of positive generators whereas strongly quasipositive knots require these conjugates to be of a special form amenable to constructing Seifert surfaces.  Positive braids are obviously strongly quasipositive (since $\sigma_{i-1,i}=\sigma_i\in B_n$). It is not obvious, yet it is true, that positive knots are strongly quasipositive (see \cite{Nakumura} or \cite{Rudolph4}).

 It is worth noting that quasipositive links are equivalent to another, more geometric class of links - the transverse $\C$-links.  These links arise as the transverse intersection of the three-sphere, $S^3 \subset \C^2$, with a complex curve.  Transverse $\C$-links include algebraic links of singularities, but are in fact a much larger class. The fact that quasipositive links can be realized as transverse $\C$-links is due to Rudolph \cite{Rudolph3} while the fact that every transverse $\C$-link is quasipositive is due to Boileau and Orekov \cite{Boileau}. 

\bigskip
 Summarizing, we have (where P stands for positive, SQP strong quasipositive, and QP quasipositive):

$$\{\mathrm{P \ braids}\}\subseteq \{\mathrm{P \ knots}\} \subseteq \{\mathrm{SQP \ knots}\} \subseteq \{\mathrm{QP \ knots}\}=\{\mathrm{transverse \ 
\C-links}\}.$$

\bigskip
\bigskip
In \cite{FourBall} \ons \ introduced an integer-valued invariant $\tau(K)$ associated to a knot, $K\subset S^3$ (see also \cite{Ras1}.)  This invariant has the following properties:

\begin{enumerate}
	\item $\tau(K_1\# K_2)=\tau(K_1)+\tau(K_2)$
	\item $\tau(\overline{K})=-\tau(K)$ with $\overline{K}$ the reflection of $K$	
	\item $|\tau(K)|\le g_4(K)$ with $g_4(K)$ the smooth slice genus of $K$
	\item $\tau(T_{p,q})=g(T_{p,q})=\frac{(p-1)(q-1)}{2}$, where $T_{p,q}$ is the $(p,q)$ torus knot and $g$ denotes Seifert genus.
\end{enumerate}

It follows that $K$ is a smooth concordance invariant.  In \cite{Livingston}, Livingston proved the following:

\begin{theorem}\label{thm:Livingston}
	(Livingston \cite{Livingston})
	Suppose $K$ is strongly quasipositive.  Then $\tau(K)=g_4(K)=g(K)$.
\end{theorem}

The statement above is slightly different than the form found in \cite{Livingston}.  Instead of requiring $K$ to be strongly quasipositive, Livingston requires that $K$ be embedded in the interior of a fiber
surface, $F$, of a torus knot.  He further requires that $K$ be null-homologous on $F$, bounding a subsurface $G\subset F$.  However, Theorem $90$ of \cite{Rudolph1} shows that these conditions are equivalent to strong quasipositivity of $K$ ($G$ is the Seifert surface for $K$ required in the definition of strong quasipositivity).

In Section $2$ we will prove:

\begin{theorem}
	\label{thm:SQP}
	Let $K$ be a fibered knot in $S^3$.  Then $\tau(K)=g_4(K)=g(K)$ if and only if $K$ is strongly quasipositive.
\end{theorem}

\begin{figure}

\begin{center}

	\psfrag{V}{$V$}
	\psfrag{K}{$K_n$}
	\psfrag{N}{$n$}

\includegraphics{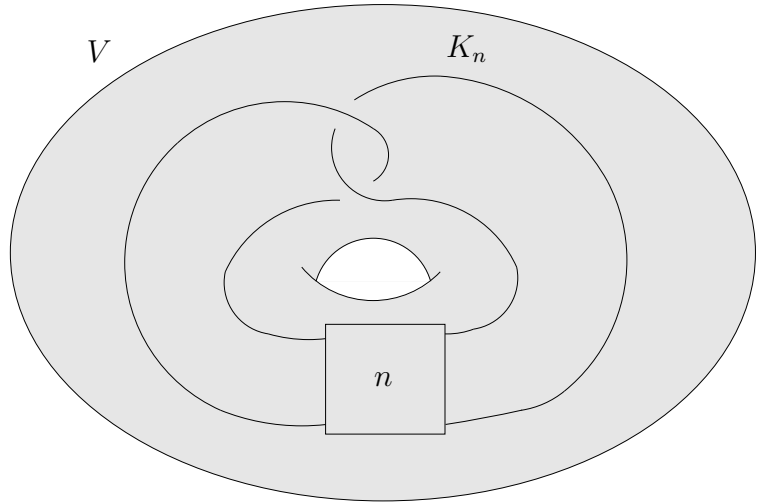}

\end{center}

\caption{\label{fig:twistknot} The twist knot $K_n$ shown in a solid torus $V$ (the box denotes $n$ full twists). $K_{-1}=$right-handed trefoil, $K_0=$unknot, $K_{1}=$figure-$8$.  The positive $n$-twisted Whitehead double, $D_+(K,n)$, of a knot $K$ is defined by identifying a tubular neighborhood of $K$ with the solid torus shown above.  $D_+(K,n)$ is the image of the twist knot $K_n \subset V$ under this identification. In the construction, one must take care to send the longitude of the solid torus above to the longitude of $K$ coming from a Seifert surface.  }

\end{figure}

It is natural to wonder if the above holds for an arbitrary knot.  We do not believe this to be the case.  Our reason follows from our ongoing study of the twisted (positive) Whitehead doubles of a knot, $K$.  Let $D_+(K,n)$ be the (positive) $n$-twisted Whitehead double of $K$ (see the caption below Figure \ref{fig:twistknot} for a definition).  Rudolph showed in \cite{Rudolph5} that  $D_+(K,n)$ is strongly quasipositive if and only if $n \le TB(K)$, where $TB(K)$ denotes the maximal Thurston-Bennequin number of $K$ \cite{EtnyreII}.  We hope to prove in an upcoming paper on the Floer homology of the Whitehead double, however, that $\tau(D_+(K,n))=g( (D_+(K,n))=1$ if and only if $n\le 2\tau(K)-1$.  In \cite{Olga}, Plamenvskaya showed that $ TB(K)\le 2\tau(K)-1$ and hence any knot for which this fails to be equality will provide a counterexample to the theorem above if the fiberedness condition on $K$ is removed (an example of such a knot is the figure-$8$, see \cite{Ng}).  For more on the \os \ concordance invariant of Whitehead doubles, see the paper of Livingston and Naik \cite{Livingston2}.

As a simple application of the theorem, we determine the strongly quasipositive iterated torus knots of type $T\{p_1,p_1 n_1+1\}\{p_2,p_2 n_2+1\}\ldots \{p_k,p_k n_k+1\}$ with $p_i\ge0$. Recall that an iterated torus knot, $T\{p_1,q_1\}\{p_2,q_2\}\ldots\{p_k,q_k\}$, is an iterated satellite knot of the unknot where the companion knot at the $i$-th stage is the $(p_i,q_i)$ torus knot.  (i.e. $T\{p,q\}$ is the $(p,q)$ torus knot, $T\{p,q\}\{r,w\}$ is the $(r,w)$ cable of the $(p,q)$ torus knot, and so on. see \cite{CableII} or \cite{Lickorish} for further explanation.) Applying the above theorem to the results of \cite{CableII} we have as corollary:

\begin{cor}
	\label{cor:Cable}
The iterated torus knot  $T\{p_1,p_1 n_1+1\}\{p_2,p_2 n_2+1\}\ldots \{p_k,p_k n_k+1\}$ is strongly quasipositive if and only if $n_i\ge0$ for all $i$.  
\end{cor}

In another direction, \ons \ have proved several remarkable theorems about the relationship between their Floer homology and the topology of lens spaces - see \cite{Lens,GenusBounds}.  For instance, \cite{Lens} places serious restrictions on the knot Floer homology invariants of any knot on which positive integral Dehn-surgery yields a lens space.  In particular, they show that any such knot satisfies $\tau(K)=g(K)$. It would follow from a conjecture of \ons \ that any such knot is fibered (see \cite{Knots}).  This conjecture is proved for an extensive list of knots which admit lens space surgeries, the so-called Berge knots \cite{Berge}. The Berge knots, in turn, are conjectured to be the only knots admitting lens space surgeries.  We have the following corollary:

\begin{cor}
	Any knot admitting a lens space surgery which appears on Berge's list is strongly quasipositive.
\end{cor}

Of course if either \ons's conjecture about the Floer homology of fibered knots or Berge's conjecture about the completeness of his list hold, then the above corollary would hold for any knot admitting a lens space surgery.  It is not true that any strongly quasipositive fibered knot admits a lens space surgery (or even an L-space surgery).  An example of such a knot is the $(2,1)$ cable of the right-handed trefoil.  It is strongly quasipositive by Corollary \ref{cor:Cable} but does not admit a lens space surgery by the examination of its knot Floer homology groups found in \cite{MyThesis}.

Regarding the relationship between quasipositivity and $\tau$, we have the following theorem. 

\begin{theorem}
	\label{thm:QP} (Plamenevskaya \cite{Olga})
		Suppose $K$ is quasipositive.  Then $\tau(K)=g_4(K)$.  
\end{theorem}

Although the above theorem does not appear in \cite{Olga}, it follows immediately from the results found there and the fact that the slice-Bennequin inequality is sharp for quasipositive knots.  This is pointed out for Rasmussen's $s$ invariant (see the remarks below) by Shumakovitch \cite{Shumakovitch}.  For completeness, we provide a proof for $\tau$ in Section $2$ and point out how this result implies that of Livingston Theorem \ref{thm:Livingston}.

Combined with Theorem \ref{thm:SQP} we then have:

\begin{cor}
	The strongly quasipositive fibered knots are precisely those quasipositive fibered knots whose smooth $4$-genera equal their Seifert genera.
\end{cor}

The question of whether this corollary holds for general knots appears to be an open problem in the subject of quasipositivity \cite{Baader}, and likely counterexamples may be the twisted Whitehead doubles of knots $K$ with $TB(K)\ne 2\tau(K)-1$ discussed above.

In addition to the corollary, Theorem \ref{thm:QP} can be used to show that a given knot is {\em not} quasipositive, for which examples seem to be lacking in the literature.  For instance we have

\noindent { \bf Examples:} The following knots are {\em not} quasipositive
\begin{enumerate}
	\item Any knot with $\tau(K)<0$ 
	\item The reflection $\overline{K}$ of any non-slice quasipositive knot, $K$
	\item Alternating knots with $\frac{\sigma(K)}{2}\ne g_4(K)$ (where $\sigma(K)$ denotes Trotter's classical signature)
	\item The twist knots $K_n$ with $n>0,n\ne 2$ (see Figure \ref{fig:twistknot})
	\item The twisted Whitehead doubles $D_+(K,n)$ of an arbitrary knot $K$ with \newline $n\ge-TB(\overline{K})$ and $n\ne b(b\pm1)$ ($n,b\in \Z$).

\end{enumerate}
		
In the above, the knots in $(2)$ (resp. $4$) are special cases of $(1)$ (resp. $3$).  
It had been previously noted by Rudolph that the twist knot $K_{1}$ (the figure-$8$) is not quasipositive \cite{Rudolph1}. There is a non-trivial overlap between $(4)$ and $(5)$ since the twist knots are twisted Whitehead doubles of the unknot, $D_+(\mathrm{Unknot},n)=K_n$. We discuss these examples at the end of the next section.

{\bf Remarks:} Though one direction of Theorem \ref{thm:SQP} follows from
Livingston's result, the proof here is independent of \cite{Livingston}.  It
should also be pointed out that Rasmussen has since discovered an invariant,
$s$, which when suitably normalized has the enumerated properties of $\tau$,
see \cite{Ras2}.  This invariant is defined using Lee's refinement \cite{Lee} of the combinatorial knot
homology theory introduced by Khovanov \cite{Khovanov}.  Rasmussen discovered
the $s$-invariant shortly after Livingston's paper appeared and indeed
Livingston's result holds for $s$ as well - his proof relies only on properties
$1-4$ of $\tau$. In contrast, Theorem \ref{thm:SQP} relies on a relationship between \os
\ theory and contact geometry which to date has not been established in
Khovanov homology.

{\bf Acknowledgement:} I would like to thank Lee Rudolph wholeheartedly for his exposition of the subject of quasipositivity and for an enlightening email discussion which illuminated a key step in the proof of Theorem \ref{thm:SQP}.  I would also like to acknowledge Charles Livingston and Olga Plamenevskaya for laying the foundations for this work in their beautiful papers \cite{Livingston,Olga}.

\section{Proof of Theorems and Corollaries}

We start this section by proving Theorem \ref{thm:SQP}.  We do this by showing a string of equivalences.  The proof relies on several major theorems and on a basic understanding of the various \os \  invariants. In particular, we exploit Giroux's theorem relating contact structures and open book decompositions and a theorem of Rudolph about Murasugi sums of quasipositive surfaces.   We do not review many of the concepts and definitions but instead refer the reader to \cite{Etnyre} for a review of contact geometry (specifically Giroux's theorem), to \cite{Survey} for an introduction to \os \ theory, and to \cite{Rudolph1} for an introduction to the different notions of positivity.

\begin{prop}
	Let $K\subset S^3$ be a fibered knot and $F$ its fiber surface. Then the following are equivalent:
	\begin{enumerate}
		
		\item $K$ is strongly quasipositive (with $F$ a quasipositive Seifert surface) 
			\item The open book decomposition associated to $(F,K)$ induces the unique tight contact structure on $S^3$
		\item  $c(\xi_K)\ne0$, where $c(\xi_K)$ is the \os \ contact invariant associated to the contact structure coming from the open book of $(F,K)$
		\item $K$ satisfies $g(K)=\tau(K).$
	\end{enumerate}
\end{prop}

\begin{proof}
	We show that each number is equivalent its predecessor.  

	\bigskip
	\noindent {$\bf 1<=>2$}  \ \  We use Giroux's fundamental theorem in contact geometry together with a theorem of Rudolph.  Associated to every fibered knot, $(Y,K)$, and its fiber surface, $F$, is an open book decomposition $(F,\phi)$ of the three-manifold $Y$.  Recall that the Murasugi sum of two surfaces along a rectangle is called plumbing.  Letting $H(+)$ (resp. $H(-)$) denote the fiber surface for the positive (resp. negative) Hopf link we say that two open books $(F_1,\phi_1)$,$(F_2,\phi_2)$ for $Y$ are stably equivalent if they become isotopic after plumbing some number of copies of $H(+)$ to each $F_i$.  That is to say,

	$$F_1\sharp k H(+) \sim F_2\sharp l H(+),$$

\noindent where $\sharp$ denotes a plumbing of surfaces and where we suppress the monodromy.  Thus open books are divided into natural equivalence classes.  Giroux's theorem states that there is a bijective correspondence between open book decompositions of $Y$ up to equivalence and contact structures on $Y$ up to isotopy.  From this point on, we will denote the contact stucture associated to a fibered knot, $K$, under Giroux's correspondence by $\xi_K$.

In the special case at hand, $Y=S^3$, foundational work of Eliashberg \cite{Eliashberg} shows that contact structures are in bijective correspondence with $\Z \cup \{pt\}$.  The exceptional point corresponds to the unique tight contact structure, $\xi_{std}$.   As for the other, overtwisted contact structures, Eliashberg shows that they are in one-to-one correspondence with homotopy classes of two-plane fields, which in turn are bijective to the integers via the Hopf invariant, $h(\xi)$.  The Hopf invariant is known to satisfy the following three properties:

\begin{enumerate}
	\item $h(\xi_{K_1\sharp K_2})=h(\xi_{K_1})+h(\xi_{K_2})$
	\item $h(\xi_{H(+)})=0$
	\item $h(\xi_{H(-)})=-1$
\end{enumerate}
		
In \cite{Rudolph2}, Rudolph shows that the Murasugi sum of two Seifert surfaces is quasipositive if and only if the two summands are quasipositive. In particular, plumbing and deplumbing with positive Hopf annuli preserves quasipositivity of the surface, and hence strong quasipositivity of the bounding knots.  Thus we see that strong quasipositivity is a characteristic of stable equivalence classes.  

We are now in a position to prove the equivalence of $1$ and $2$ stated in the proposition.  Suppose first that we are given a knot $K$ such that $\xi_K=\xi_{std}$.  Since the unknot induces $\xi_{std}$, Giroux's theorem says that $K$ is stably equivalent to the unknot and Rudolph's theorem in turn shows that $K$ is strongly quasipositive.

Now assume conversely that $\xi_K\ne\xi_{std}$ (i.e. that $\xi_K$ is an overtwisted contact structure).  By the above remarks, we must only exhibit a knot $K'$ which is stably equivalent to  $K$ and which is not strongly quasipositive.  Since the overtwisted contact structures are distinguished by $h(\xi)$, and by the additivity of $h(\xi)$ (Property 1), it suffices to show that the stable equivalence class of the overtwisted contact structure with $h(\xi)=0$ contains a non-strongly quasipositive representative.  To this end, take any representative $L$ for the overtwisted contact structure with $h(\xi_L)=1$.  Now plumb one negative Hopf link to $L$, together with some positive Hopf links (to ensure that the result is a knot) to obtain,

$$K'= L\sharp H(-) \sharp k H(+),$$

\noindent so that

$$h(\xi_{K'})= h(\xi_L)+ h(\xi_{ H(-)}) + k h(\xi_{H(+)}) = 1-1+k\cm0 =0.$$

	\noindent Using Rudolph's theorem again, we see that $K'$ is not strongly quasipositive, since the negative Hopf link is not.  This proves the equivalence of $1$ and $2$.
	
\bigskip
	
\noindent {\bf Remark:} We are indebted to Lee Rudolph for suggesting the above idea.

	\bigskip
	\noindent	{$\bf 2<=>3$} \ \ To a closed oriented three-manifold, $Y$, and Spin$^c$ structure, $\spinc$,  \ons \  introduced a chain complex $\CFa(Y,\spinc)$ \cite{HolDisk}.  A null-homologous knot $(Y,K)$ induces a filtration $\Filt(Y,K,i)$ of this chain complex, i.e. there is an increasing sequence of subcomplexes:

$$ 0=\Filt(Y,K,-i) \subseteq \Filt(Y,K,-i+1)\subseteq \ldots \subseteq
\Filt(Y,K,n)=\CFa(Y,\spinc).$$
	
\noindent For the definition of this filtration see \cite{Knots,Ras1}.  In \cite{Contact}, \ons \  showed that for fibered knots, $H_*(\Filt(-Y,K,-g))\cong \Z$, where $-Y$ indicates the manifold $Y$ with opposite orientation, and where $g$ indicates the Seifert genus of $K$.  They further showed that $H_*(\Filt(-Y,K,i))\cong 0$ for $i<-g$. Let $c_0(K)$ denote a generator for $H_*(\Filt(-Y,K,-g))$, and let $c(K)$ denote its image under the map on homology induced from the natural inclusion map

$$i : \Filt(-Y,K,-g) \rightarrow \CFa(-Y).$$

Theorem 1.3 of \cite{Contact} states that if the contact structures induced by the open books of two fibered knots $K_1$ and $K_2$ are equivalent, then the invariants $c(K_1)$ and $c(K_2)$ are the same, up to sign.  It follows that $c(K)\in \HFa(Y)/(\pm1)$ is an invariant of the contact structure associated to the open book for $K$, and we denote this contact invariant by $c(\xi_K)$.  Theorem 1.4 of \cite{Contact} states that this invariant vanishes if the contact structure $\xi_K$ is overtwisted.   It is not difficult to see that the invariant associated to $(S^3,\xi_{std})$ is equal to a generator for $\HFa(S^3)\cong \Z$.  Any other contact structure on $S^3$, being overtwisted, must necessarily have zero invariant by \ons's theorem.

\bigskip 
	\noindent {$\bf 3<=>4$} For a knot $K\subset S^3$, the invariant $\tau(K)$ is defined as:

$$\tau(K)=\mathrm{min}\{j\in\Z|i_* : H_*(\Filt(S^3,K,j))\longrightarrow \HFa(S^3)\cong \Z \  \mathrm{is \ non\-trivial}\}.$$

The above discussion shows that contact invariant of $K$ is non-zero if and only if $i_* : H_*(\Filt(-S^3,K,-g) \rightarrow \HFa(-S^3)$ is non-trivial.  Now the image of $K$ under an orientation reversing diffeomorphism of $S^3$, is the reflection, $\overline{K}$.  Thus, we have

$$\tau(\overline{K})=\mathrm{min}\{j\in\Z|i_* : H_*(\Filt(-S^3,K,j))\longrightarrow \HFa(-S^3)\cong \Z \  \mathrm{is \ non\-trivial}\},$$

\noindent and we see that the contact invariant of $\xi_K$ is non-zero if and only if $\tau(\overline{K})=-g$.  But Property $(2)$ of $\tau$ stated in the introduction implies that $\tau(K)=-\tau(\overline{K})=g$.  This completes the proof.  

\end{proof}

\noindent {\bf Proof of Corollary \ref{cor:Cable}:} The result follows immediately from Theorem 1.3 of \cite{CableII} which states that an iterated torus knot $K=T\{p_1,p_1 n_1+1\}\{p_2,p_2 n_2+1\}\ldots \{p_k,p_k n_k+1\}$ satisfies $\tau(K)=g(K)$ if and only if  $n_i\ge0$ for all $i$ (assuming positive $p_i$). 

\bigskip

We now turn our attention to Theorem \ref{thm:QP}

\noindent {\bf Proof of Theorem \ref{thm:QP}:} The celebrated slice-Bennequin inequality proved by Rudolph \cite{Rudolph6,Rudolph5} states that, for a Legendrian knot $K\subset (S^3,\xi_{std})$

 $$tb(K) + |rot(K)| \le 2g_4(K)-1, $$

\noindent Where $tb(K)$ and $rot(K)$ denote the Thurston-Bennequin and rotation numbers of $K$, respectively (see \cite{Etnyre} for definitions and an introduction to Legendrian knots).  On the other hand, Plamenevskaya proves in \cite{Olga} that

$$ tb(K) + |rot(K)| \le 2\tau(K)-1. $$

\begin{figure}

\begin{center}

\includegraphics{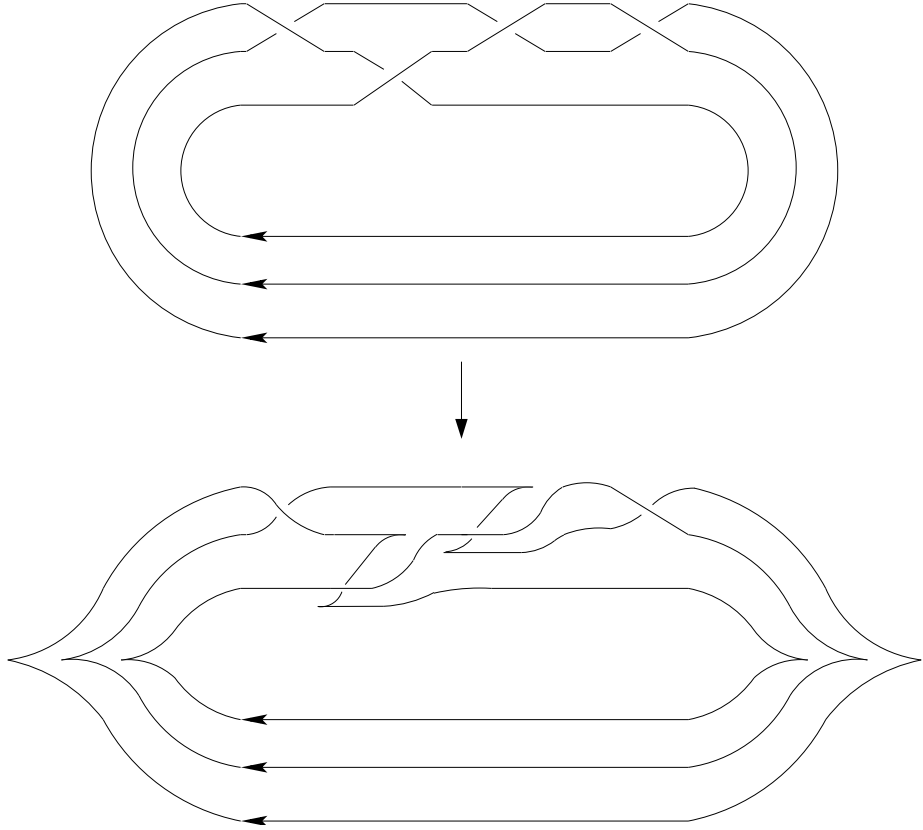}

\end{center}

\caption{\label{fig:Legendrianize} Depiction of the ``Legendrianization'' of the braid $\sigma_1\sigma_2^{-1}\sigma_1^{-1}\sigma_{1}$.  The cusps on the left hand side of the page are so-called ``up left'' cusps with respect to the orientation on $K$ shown.}

\end{figure}

\noindent Since $|\tau(K)| \le g_4(K)$, the theorem will follow if we can show that the slice-Bennequin inequality is sharp for quasipositive knots. (Indeed, this is apparent from \cite{Rudolph6,Rudolph5} but we prove it here for the reader's convenience).  We prove this by explicit construction.  Figure \ref{fig:Legendrianize} depicts an algorithm which, given a knot, $K$, presented as a braid closure, constructs a Legendrian representative, $\tilde{K}$.  Now the well-known formula for $tb(K)$ states

$$ tb(K) = \mathrm{writhe}(K)-\#\{\mathrm{left \ cusps}\}, $$

\noindent whereas the absolute value of the rotation number is given by

$$ |rot(K)|= |\#\{\mathrm{down \ left \ cusps}\}-\#\{\mathrm{up \ right \ cusps}\}|, $$

\noindent where a cusp is said to be a down (resp. up) cusp if the direction of travel along the knot (with respect to a given orientation of K) is down (resp. up) as we traverse the cusp.

From Figure \ref{fig:Legendrianize} it is apparent that a positive generator $\sigma\in B_n$ adds $1$ to $tb(\tilde{K})$ and does not affect $|rot(\tilde{K})|$, whereas a negative generator $\sigma^{-1}$ adds $-2$ to $tb(\tilde{K})$ and adds $1$ to $|rot(\tilde{K})|$.  Let $n_{+}$ (resp. $n_{-}$) denote the number of positive (resp. negative) generators in a given braid $\beta$, and let $b$ denote the braid index.  If the braid happens to be quasipositive, $\beta = \Pi_{k=1}^m w_k \sigma_{i_k} w_k^{-1}$, then $n_{+}=n_{-}+m$.  Thus

$$ tb(\tilde{K}) + |rot(\tilde{K})| = -b + n_+ - 2n_- + n_- = -b +n_+ - n_- = -b + m.$$

Rudolph \cite{Rudolph3}, however, constructs a smooth complex curve which intersects $S^3$ transversely in the closure of the given quasipositive $\beta$, whose genus is $\frac{-b+m+1}{2}$.  Thus we have shown that the slice-Bennequin inequality is sharp for quasipositive knots.  In the special case that $\beta$ is strongly quasipositive, the genus of the quasipositive Seifert surface which the knot bounds (shown in Figure \ref{fig:SS}) is also  $\frac{-b+m+1}{2}$, and hence we recover Livingston's result, Theorem \ref{thm:Livingston}.

\bigskip

We conclude by explaining the examples at the end of the introduction.  The fact that a knot with $\tau(K)<0$ is not quasipositive follows immediately from Theorem \ref{thm:QP}, as $g_4(K)\ge0$.  Since Theorem \ref{thm:QP} shows that any non-slice  quasipositive knot satisfies $\tau(K)=g_4(K)>0$, we see that its reflection must satisfy $\tau(\overline{K})=-\tau(K)<0$, and hence cannot be quasipositive.   \ons \ prove in \cite{AltKnots} that $\tau(K)=\frac{\sigma(K)}{2}$ for any knot which admits an alternating projection.  Example $3$ of the introduction follows immediately.  The twist knots are alternating and satisfy $\sigma(K_n)=0$ for all $n>0$. It follows from the work of Casson and Gordon \cite{Casson} that the only one which is slice is $K_{2}$ ($6_2$ in the tables), and so we obtain Example $4$.  As for the Whitehead doubles, Livingston and Naik \cite{Livingston2} show that $\tau(D_+(K,n)=0$ if $n\ge -TB(\overline{K})$.  When, in addition, $n\ne b(b\pm1)$, the Alexander polynomial,

$$\Delta_{D_+(K,n)}(T)=-nT+(2n+1)-nT^{-1},$$

\noindent is not of the form, $F(T)F(T^{-1})$, which is a necessary condition for a knot to be slice \cite{Lickorish}.

\commentable{
\bibliographystyle{plain}
\bibliography{biblio}
}

\end{document}